\input psfig.sty 

\hsize=6truein
\overfullrule=0pt

\font\title=cmr10 scaled\magstep4
\font\subtitle=cmr10 scaled\magstep2
\font\author =cmr10 scaled\magstep1

\font\eaddr=cmtt10
\font\ninerm=cmr8
\baselineskip=16pt plus 1pt minus 1pt
\def\QED{\vrule height8pt width4pt depth0pt}

\def\proclaim #1: #2\par{\medbreak
	      \noindent{\bf#1:\enspace}{\sl#2}\par\medbreak}
\pageno=1	      

\hfill March 18, 2013\break

\vskip 1truein
\centerline{\title Sequences with long range exclusions}
\vskip .7truein
\centerline{{\author Kari Eloranta}\footnote{$^*$}
{Research partially supported by the Magnus Ehrnrooth foundation}}
\vskip .1truein
\centerline{\author Institute of Mathematics}
\centerline{\author Aalto university}
\centerline{\author Finland}
\vskip .1truein
\centerline{\eaddr kari.eloranta@aalto.fi}


\vskip .3truein
\centerline{\subtitle Abstract}
\vskip .3truein
\centerline{\vbox{\hsize 4.5in \noindent \ninerm \strut Given an alphabet $S$, we consider the size of the subsets of the full sequence space $S^{\rm {\bf Z}}$ determined by the additional restriction that $x_i\not=x_{i+f(n)},\ i\in {\rm {\bf Z}},\ n\in {\rm {\bf N}}.$ Here $f$ is a positive, strictly increasing function. We review an other, graph theoretic, formulation and then the known results covering various combinations of $f$ and the alphabet size. In the second part of the paper we turn to the fine structure of the allowed sequences in the particular case where $f$ is a suitable polynomial. The generation of sequences leads naturally to consider the problem of their maximal length, which turns out highly random asymptotically in the alphabet size.}}

\vskip .4truein
\centerline{\vbox{\hsize 4.5in \noindent Keywords: symbolic dynamics, chromatic number, formal language, long range order, difference sets, Poincar\'e recurrence, Pascal's simplex, compound geometric distribution. \hfill\break \hfill\break
2010 AMS Classification: 37B20, 05C15, 11B05, 60F05 \hfill\break
Running head: Sequences with long range exclusions}}

\vfill
\eject

$ $
\vskip .4truein
\noindent {\subtitle 0. Introduction}
\vskip .2truein 

\noindent Given an {\bf alphabet} $S=\{1,2,3,\ldots, d\}$ on $d$ symbols the {\bf full shift} is the dynamical system $\left(S^{\rm {\bf Z}}, \sigma\right).$ Here $\sigma$ is the left shift: $\sigma(x)_i=x_{i+1},\ \forall i\in {\rm {\bf Z}}.$ As the full shifts are rather simple objects the main study in topological dynamics usually concentrates on subshifts of finite type (SoFT), where one further restricts the sequences by listing all allowed neighboring symbols. This leads to the transfer matrix formulation which for 1-dimensional sequences has a pretty complete theory.

If the symbol interaction is not nearest neighbor but finite range, by using a larger alphabet (from blocks of $S$-symbols) one can still remain in the SoFT setup. For infinite range interactions this is not anymore possible and the transfer matrix formalism breaks down.

Suppose that we impose an infinite restriction, an {\bf exclusion rule}, and define:
$$X_{(d,f)}=\left\{x\in S^{\rm {\bf Z}}|\ x_i\not=x_{i+f(n)},\ i\in {\rm {\bf Z}},\ n\in {\rm {\bf N}}\right\}, \leqno{(0.1)}$$
where $f:{\rm {\bf N}}\rightarrow {\rm {\bf N}},\ ({\rm {\bf N}}=\{1,2,3,\ldots\})$ is strictly increasing hence also unbounded. Note that the definition implies a symmetric rule (w.r. to any $i$). One could interpret this as some kind of exclusion principle: symbols on \lq\lq shells\rq\rq\ of radius $f$ repel their kind. Mike Keane proposed a model like this, with the particular choice $f(n)=n^2$ ([LM]). As a first question he asked for which $S$ the space is nonempty. Turns out that even for the simplest choices of the {\bf jump sequence} $f$ answering this is fairly involved. Note that in symbolic dynamics context it would be most natural to expect a finite alphabet, $d<\infty.$

Alternatively one could view ${\rm {\bf Z}}$ {\bf as a graph} where there is an edge between two distinct vertices/integers $i$ and $j$ iff $|i-j|=f(n)$ for some $n\in {\rm {\bf N}}.$ Coloring the graph by requiring that adjacent vertices have different colors leads to the natural question of the {\bf chromatic number} of this (undirected) graph. This formulation seems to be due to Paul Erd\"os who was curious about the $f$ for which the chromatic number is finite.

In this paper we will first present some examples and then review the known results. Through these the problem is then connected to more general dynamics and to additive number theory. It turns out that in some sense for most combinations of $d$ and $f$ the shift spaces are empty. This will lead to us to consider how the termination of the sequences comes about. Analysis of these generation algorithms is the main content of the latter half of the paper.

\vskip .5truein
\noindent {\subtitle 1. Infinite sequences}

\vskip .2truein
\noindent {\subtitle 1.1. First observations}
\vskip .2truein

\noindent Together with the two-sided sequences of (0.1) we also deal with the one-sided version $X_{(d,f)}^+$ where the rule is restricted to $S^{\rm {\bf N}}.$ If we can show that $X_{(d,f)}^+$ is small/empty then by $X_{(d,f)}\big|_{{\rm {\bf N}}}\subset X_{(d,f)}^+$ (the restriction to positive integers) the same holds for the two-sided space. The inclusion holds since in $X_{(d,f)}^+$ the past haunts only back to the coordinate 0 whereas in $X_{(d,f)}$ the block can of course arise from a symbol arbitrarily far in the past. A second observation is even more basic but we state it nevertheless:  $X_{(d,f)}\subset   X_{(d'f)}$ if $d\le d'.$ This holds for one-sided sequences, too. If we can show the smallness for some alphabet size, it holds for all smaller alphabets.

Let us first consider some explicit cases to gain some insight into the contribution of the alphabet size and growth rate of the jump sequence. By {\bf lexicographically generating} the one-sided sequences we mean that one proceeds from $x_1=1$ rightwards to the nearest coordinate always assigning the smallest allowed symbol. 

\vskip .2 truein
\noindent {\bf Example 1.1. (linear $f$):} Consider the case $S=\{1,2\}$ and $f(n)=2n.$ Clearly if $x_0=1$ then $x_{2k}=2, \forall k\not=0$ but $x_2=2$ implies $x_{2m}=1, \forall m\not=1,$ a contradiction. So $X_{(2,2n)}=\emptyset.$

On the other hand $X_{(2,2n-1)}\not=\emptyset$: it consists of the two periodic points $(12)^*.$

It is worth noticing that the emptiness is not due to too small $S.$ Let $f(n)=kn,\ k\ge 2$, fixed integer. If $x_0=1$, then $x_{kn}=\neg 1, \forall n\ge 1$ ($\neg s$ meaning \lq\lq any symbol but $s$\rq\rq). If then $x_k=2$, we have $x_{lk}\notin \{1,2\},\ \forall l\ge 2.$ Continuing like this considering symbol choices at coordinates that are multiples of $k$ one exhausts any finite alphabet i.e. $X_{(d,kn)}=\emptyset$ for any finite $d.$ 

\vskip .2 truein
\noindent {\bf Example 1.2. ($d=2$ and power $f$):} Suppose $S=\{1,2\}$ and $f(n)=n^r,\ r=2,3,\ldots$ Let $x_0=1$ then $x_{2i}=1,\ \forall i\in {\rm {\bf Z}}$ so in particular $x_{2^r}=1$ which generates a contradiction. Hence $X_{(2,n^r)}=\emptyset.$ Here the minimal size of $S$ is crucial: it strongly forces parity.

\vskip .2 truein
\noindent {\bf Example 1.3. (fast growing $f$):} Suppose there is a natural number $m$ which does not divide $f(n)$ for any $n=1,2,\ldots$ (subsequently we use the notation $a|b$ meaning $a$ divides $b$). Clearly for e.g. $n^r,\ n^n, n!\ldots$ this isn't the case, but $p^n$ has this property for prime $p.$ For $d\ge m$ one can then have infinite periodic points. Hence for example $X_{(3,2^n)}\not=\emptyset$ (at least 6 elements, $(123)^*$ and its permutations).

It is straightforward to see that $X_{(d,n!)}^+=\emptyset$ at least for $d=2,3.$ One simply sets $x_1=1,$ lays down the blocks from it, sets $x_1=2$ and then the blocks that it generates etc. With $d=2$ one is stuck in the second assignment round, with $d=3$ in the sixth.

 $X_{(4,n!)}^+$ is already considerably more interesting. With computer one can see that lexicographically generated, the sequence quickly becomes periodic with period $(1234123412312341231234123)$ (of length 25). However this period cannot persist indefinitely since starting with $x_1=1$ we would have $x_{7!}=1.$ The sequence does not terminate then though but instead there is a short transient after which the original period is reinstated (in the lexicographic generation of the sequence). Since $25\big|n!,\ \forall n\ge 10$ we see that this periodicity would be interrupted infinitely often. Is there an infinite sequence in  $X_{(4,n!)}^+$? In computer runs one with the given period will persist at least for 500.000 steps.

\vskip .2truein
\noindent The existence or non-existence of a period leads naturally to a language theoretic characterization of these spaces. For a background on these we refer to [HU]. Let $|s|$ denote the length of the string/sequence $s$. 

\proclaim Theorem 1.1.: If for any $m\in  {\rm {\bf N}}$ there is $n\in  {\rm {\bf N}}$ such that we have $m\big|f(n)$ then the words satisfying the exclusion rule of (0.1) do not form a context-free language. Hence the strings do not form a regular language (sofic shift) either. 
\par

\noindent {\bf Proof:} A context-free language necessarily satisfies a Pumping lemma: any sufficiently long string $s$, say $|s|\ge k$, can be written as $s=uvxyz$ such that (i) $|vxy|\le k$, (ii) $|vy|\ge 1$ and (iii) $uv^nxy^nz$ is an allowed string for all $n\in {\rm {\bf N}}.$ Here (ii) means that the \lq\lq pumping\rq\rq\ is indeed nontrivial action. If $\min\{|v|,|y|\}=0$ but (ii) still holds, this more general Pumping lemma for context free languages is reduced to that for regular languages.

In our context the validity of a Pumping lemma would mean that there exists some period string that allows arbitrarily long repetitions always resulting in a legal string. But if this period is of length $m$ and $m\big|f(n)$ for some $n\in {\rm {\bf N}}$ then e.g. the first symbol in the block contradicts the identical symbol in the $(f(n)/m)^{\rm th}$ block. \hfill\QED  

\vskip .2truein
\noindent {\bf Remarks:} The result ties the divisibility properties of of $\left\{f(n)\right\}_{n\in {\rm {\bf N}}}$ to the complexity of the sequences: if non-divisibility prevails periodic sequences exist, if not then the recipe for constructing infinite sequences must be more complicated. Examples like $X_{(3,2^n)}$ or $X_{(4,{\rm primes})}^+$ containing periodic points are rare specialities. But the latter half of Example 1.3. suggests that forbidden periodicity is likely just part of the story.

\vskip .5truein
\noindent {\subtitle 1.2. Lacunary case}
\vskip .2 truein

\noindent Considering one-sized sequences it is clear that for a finite $d$ a faster growing $f$ allows more freedom for the sequence generation. In particular P. Erd\"os seems to have formulated in the mid 80's the question whether lacunary sequences are sufficiently fast growing for infinite sequences/finite chromatic number to prevail. Recall that $\{x_n\}$ is {\bf lacunary} if there is $\epsilon>0$ such that $x_{n+1}/x_n\ge 1+\epsilon,\ \forall n.$ This problem was solved by Y. Katznelson and published later in [K]:

\proclaim Theorem 1.2.: For any lacunary $f$ there is a finite $d$ such that $X_{(d,f)}\not=\emptyset.$\par  

\noindent {\bf Remarks 1.} This is the sequence translation of the original solution which was in terms the chromatic number. The solution connected the problem to another one, in Diophantine approximation (also by Erd\"os), and solved them using a novel dynamcial system formulation.\hfill\break
\noindent {\bf 2.} Katznelson didn't provide an explicit bound for the chromatic number/number of symbols needed although one was implicit in his work. Other have later improved the bound the best being now $1+c|\log\epsilon|/\epsilon$ where $\epsilon$ is as in the definition of lacunarity and $c$ an absolute constant ([PS]). Note that while of course useful this still doesn't resolve explicit cases like Example 1.3. i.e. whether $X_{(4,n!)}^+$ is non-empty.  

\vskip .2truein
\noindent It is natural to connect this to recurrence as follows (see [W] for more detail). Given $(d,f)$ let $X$ be the subset (0.1) and $\rho$ the natural metric it inherits from the full shift.

\proclaim Definition 1.3.: $D\subset {\rm {\bf N}}$ is a {\bf non-recurrence set} for the (topological) dynamical system $(X,\sigma,\rho)$ if for some positive $b$ and all $x\in X$ and $d\in D$ it holds that $\rho\left(\sigma^d x,x\right)\ge b.$ \par

\noindent It is easy to show that the finiteness of the chromatic number is equivalent to $D$ being a non-recurrence set ([W]). So in particular Theorem 1.2. shows that a lacunary $D$ is a set of non-recurrence. Not much seems to be currently known about non-lacunary non-recurrence sets - very few examples are known.

For subsequent use will be to define $D$ a {\bf set of recurrence} if $X$ has infinite chromatic number.

\vskip .5truein
\noindent {\subtitle 1.3. Sublacunary realm}
\vskip .2 truein
\noindent {\subtitle 1.3.1. Dynamics approach}
\vskip .2 truein

\noindent We now chart a bit what is known about the sublacunary case. The resulting map will unfortunately be somewhat patchy, perhaps indication of the fact that not all right concepts are yet known. Clearly recurrence sets and non-recurrence sets partition the space of sequences but unfortunately there doesn't seem to be a definitive criterion for the recurrence sets. Within them most attention has been paid to Poincar\'e sequences.

\proclaim Definition 1.4.: A sequence $D\in {\rm {\bf N}}$ is a {\bf Poincar\'e sequence} if for a (measurable) dynamical system $(Y,{\cal B},\mu, T)$ and any $B\in {\cal B},\ \mu(B)>0$ it holds that $\mu\left(T^{-d_k}B\cap B\right)>0$ for some $d_k\in D.$\par

\noindent {\bf Remarks: 1.} The Poincare Theorem concerns the case $D={\rm {\bf N}}$, showing this $D$ to be a set of recurrence.\hfill\break
{\bf 2.} An ingenious example by I. Kriz shows that Poincare sequences are not all of recurrence sets.

\vskip .2truein
\noindent  As the Definition above is not explicit, other forms have been searched for. For a subset of integers $A$, let the {\bf difference set} be $A-A=\{a-a'|\ a,a'\in A\}$ and let $A^{(N)}=A\cap\{1,2,3,\ldots,N\}.$ By the {\bf upper density of} $A$ we mean $\limsup_{N\rightarrow\infty}{|A^{(N)}|\over N}.$ Betrand-Mathis established the equivalence: $D$ is a Poincar\'e sequence if and only if 
$D\cap (A-A)\not=\emptyset$ for all $A$ of positive upper density ([B-M]). Due to this formulation, in the literature such $D$ is also called an {\bf intersective sequence}.

\vskip .2truein
\noindent The developments of the subject gained definitive momentum when Lazlo Lov\'asz conjectured (apparently just verbally) that if one insists that $n^2\notin A-A$ for all natural numbers $n$ then $A$ cannot have positive density. In the late seventies Furstenberg and S\'ark\"ozy separately managed to prove (in [F] and [S]):

\proclaim Theorem 1.5.: Given $\delta>0$ there is $N_0(\delta)$ such that if $N\ge N_0(\delta)$ and $|A^{(N)}|\ge \delta N$ then there is natural $n$ such that $n^2\in A-A.$
\par

\noindent {\bf Remark:} The proofs were of different character; Furstenberg's ergodic theoretic, S\'ark\"ozy's utilized Fourier analytic techniques enabling quantitative estimates for the density. This will have some implications in our further work.

\vskip .2truein
\noindent Theorem 1.5. together with the Betrand-Mathis -equivalence shows that $\{n^2\}$ is an example of a Poincar\'e/intersective sequence. Later on more examples of such sequences have been worked out by Furstenberg and many others. We are not going to list these but just highlight that e.g. all monomials $n^k$ and more generally polynomials satisfying the hypothesis of Theorem 1.1. are intersective (the latter result is due to Kamae and Mend\'es France, [KMF]).

While Betrand-Mathis alternative is conceptually valuable, a criterion often more usable is

\proclaim Lemma 1.6.: If for all $\alpha\in (0,2\pi)$ $\lim_{n\rightarrow\infty}{1\over n}\sum_{k=1}^ne^{i\alpha s_k}$ then $\{s_k\}$ is a Poincar\'e sequence.\par

\noindent {\bf Remark:} These are called {\bf Weyl sums} due to his first use of them to show uniform distribution (for integers).

\vskip .2truein
\noindent The results above have direct bearing on our original problem. Consider for simplicity its 1-sided version with $d$ symbols and $f$ an intersective polynomial. Let $A_i=\{j|\ x_j=i\},$ these are just the sets of coordinates where the $i^{\rm th}$ symbol is found. Suppose now that  $A_i,\ i=1,\ldots,d$ partition ${\rm {\bf N}}$ i.e. the 1-sided sequence is uniquely defined everywhere on ${\rm {\bf N}}$. Then Theorem 1.5. generalized as indicated above implies that if the exclusion is to hold for this sequence, for sufficiently large $N$ the upper densities of $A_i$'s cannot add up to 1, a contradiction. Hence

\proclaim Theorem 1.7.: For intersective polynomials $f$ the spaces $X_{(d,f)}^+$ and therefore all $X_{(d,f)}$ are empty. \par

\vskip.2truein
\noindent To gain some insight on what can be done from first principles and more importantly, how the sequence termination takes place, let us take for a moment a purely elementary view. 

\vskip .2truein
\noindent {\bf Example 1.4. ($d=3,4$ and $n^2$):} To directly resolve the case $d=3$ one can utilize simple identities like $3^2+4^2=5^2$: instead of dealing with all the entries upto a given distance left and right of origin, one only needs to check a rather thin subset of them. Let $x_0=1$ and generate its blocks on $\{-25,\ldots,\ 25\}.$ Then set $x_{16}=2$ and generate the blocks from it. In third iteration one sets $x_{-9}=3$ (or $x_{25}=3$ just as well) since there is a double block at it already due to the given identity. Now the full alphabet is in use. If e.g. at the doubly blocked coordinate 7 one sets $x_7=1$ the following (fourth) assignment iterate succeeds but next one results in a full block at several sites on $\{-25,\ldots,\ 25\}.$ Hence $X_{(3,n^2)}=\emptyset.$

This procedure of trying to generate a contradiction in the construction of a two-sided sequences gets rather unwieldy for larger $S$ even for this $f.$ Suspecting the spaces might be empty one could switch instead to considering one-sided sequences i.e. aiming to show that $X_{(4,n^2)}^+=\emptyset.$ Indeed here one succeeds with a {\bf computer assisted proof} (CAP). Checking all these sequences systematically is a manageable task and turns out they are finite, the longest being of length 47. $X_{(5,n^2)}^+.$ is exponentially harder and CAP on a desktop machine does not seem to terminate (sequences of length about 170 exist).

The problem gets quickly rather intractable for higher monomials $n^r,\ r\ge 3$ one reason being that there are no simplifying identities like $k^r=n^r+m^r$ available then due to Fermat's theorem. For $r=3$ a CAP program finds one-sided sequences of length at least 300. Note that by the argument in Theorem 1.1. we immediately know that there are no periodic sequences in any of the spaces $X_{(d,n^r)}.$

\vskip .5truein
\noindent {\subtitle 1.3.2. Additive combinatorics}
\vskip .2truein

\noindent S\'ark\"ozy's result inspired plenty of subsequent additive number theory. Some of these results shed light to our original problem. They were often parallel results to those in the previous section, while simultaneously being quantitative. We will now briefly highlight this development.

\vskip .2truein
\noindent S\'ark\"ozy's proof of Theorem 1.5. utilized the Hardy - Littlewood  Circle Method. By honing it further Balog, Pelik\'an, Pintz and Szemer\'edi established in 1994

\proclaim Theorem 1.8.: Fix a natural $k\ge 2.$ If $n^k\notin A-A$ for all $n$ then
$${|A^{(N)}|\over N} \ll {1\over {(\log N)^{c\log\log\log\log N}}}.$$
\par

\noindent  Here $\ll$ means \lq\lq less than constant times\rq\rq.  Best $c$ has been worked out to be $1/\log{3}.$ This in turn has been later been extended to {\bf intersective} polynomials $f\in {\rm {\bf Z}}[x]$ perhaps culminating to J. Lucier's result in 2006 ([L]):

\proclaim Theorem 1.9.: Suppose $f\in {\rm {\bf Z}}[x]$ is an intersective polynomial of degree $k\ge 2$ with positive leading term. If $f(n)\notin A-A$ for all $n$ with $f(n)>0$ then
$${|A^{(N)}|\over N} \ll \left({{(\log\log N)^\mu}\over{\log N}}\right)^{1\over {k-1}}$$
for $\mu=3$ if $k=2$, $\mu=2$ if $k>2$ and the constant only depends on $f.$\par

\vskip .2truein
\noindent Like in the dynamics approach, from these results one can readily deduce Theorem 1.7. Additionally here one gets a upper bound for the instant when the \lq\lq total density\rq\rq\ of the $d$ subsequences drops below 1 i.e. the contradiction takes place. But as one can see from the formulas the bound will be huge. This is due to the method of proof, not the phenomenon at hand. This absence of reasonable bounds will motivate us to investigate in the next section how long sequences satisfying (0.1) in the sublacunary realm can actually prevail. 

\vskip .2truein
\noindent One could of course argue that in the sublacunary case there is nothing special about in polynomials. Perhaps so, in which case the findings above are just to indicate explicit cases, where the original problem can be solved. However there is a reason to believe that Theorem 1.7 is generic in the sense that for all sublacunary rates of growth, there are $f$ such that the chromatic number is infinite i.e. $X_{(d,f)}$ is empty for all $d.$ (see [AHK]).

Our knowledge of the sublacunary realm is hence left in a somewhat curious state. We do not know whether for all sublacunary growth rates there are non-empty $X_{(d,f)}$ coexisting with the empty ones. Our Example 1.1. ($\{2n\}$ is recurrent, $\{2n-1\}$ in non-recurrent) shows that this is possible for at least some rate. Is there an upper bound for the growth rate below which such coexistence is possible? Perhaps the inevitable conclusion from this state of affairs is that the growth rate is just too crude a measure for the problem. One should use something (explicit!) that genuinely captures the recurrence-properties of the $f$-sequence.

\vskip .5truein
\noindent {\subtitle 2. Finite sequences}
\vskip .2truein

\noindent {\subtitle 2.1. Random generation}
\vskip .2truein

\noindent We now proceed to the random generation of the simpler one-sided subset of sequences compatible with the exclusion (0.1) and the analysis of the outcome of this algorithm. The main results here are aimed at showing that a probabilistic model gives surprisingly accurate information on the length distribution. While the methods are applicable to a wide variety of jump sequences $f=n^2$ is the only one that we give a statistical assessment.

\vskip .2truein
\noindent Our algorithm generates random sequences obeying the restriction (0.1) on $S^{\rm {\bf N}}.$ Each coordinate is chosen independently and uniformly but in such a way as to respect the restrictions from all the coordinates on its left. 

\vskip .2truein
\noindent {\bf Algorithm v2.0:}\hfill\break 
\indent {\tt 0. set} $M\ge 1,$ {\tt let} $S_j=S$ {\tt at each}  $j\in \{1,\ldots,M\}$ {\tt and set} $i=1.$\hfill\break
\indent {\tt 1. if $S_i=\emptyset$ then {\bf halt},\hfill\break
\indent\hskip.75cm else pick uniformly a random symbol} $s\in S_i.$ \hfill\break
\indent {\tt 2. update} $S_j\leftarrow S_j\setminus \{s\}$ {\tt for\ all} $j=i+f(n)\in \{i+1,\ldots,M\},\ n\in {\rm {\bf N}}.$\hfill\break
\indent {\tt 3. if } $i=M$ {\tt halt\ and\ call} {\bf full\ length},\hfill\break
\indent\hskip.75cm {\tt else} $i\leftarrow i+1$ {\tt and go to 1.}

\vskip .2truein
\noindent So just halting means that the random procedure did not result in a legal string of symbols of length $M$ but produced a contradiction earlier. \lq\lq Full length\rq\rq\ means that length $M$ was reached and the sequence needs to be continued with larger $M$ to further address its viability.

One might perhaps also want to think the the sequence $\{s_i\}$ as a {\sl random walk in a random environment (RWRE)}. Since the walker itself generates the future obstacles one could call the walk in this sense {\sl self-similar}.  
    
\vskip.2truein
\proclaim Definition 2.1.: Given a site $j\ge 2$ call the set of sites $D_j=\{1,2,\ldots,j-1\}\cap\{k\ |\ k=j-f(n)\ge 1,\ n\in {\rm {\bf N}}\}$ the {\bf dragnet} of $j$ and let $D_1=\emptyset$. The cardinality of the dragnet is an increasing step function. Let the interval between the $(d+i-1)^{\rm st}$ and $(d+i)^{\rm st}$ steps span the $i^{\rm th}$ {\bf interval}, $i\ge 1.$ It's length is denoted by $l_i.$\hfill\break
If the sites in the dragnet $D_j$ support the entire alphabet $S$, the site $j$ has {\bf full block}. The site where the first full block can materialize is the {\bf first jump}, the next step after that is {\bf second jump} etc. \par

\vskip .2truein
\noindent Suppose that we take the particular choice $f(n)=n^2.$ Then we have explicitly
 
\proclaim Proposition 2.2.: Given an alphabet of size $d$, the location of the $i^{\rm th}$ jump is at $(d+i)^2+1,\ i=1,2,\ldots$ The $i^{\rm th}$ interval $\{(d+i-1)^2+1,\ldots,(d+i)^2\}$ is of length $l_i=2(d+i)-1.$ \par

\vskip .1truein
\noindent The proof is simple calculation and we omit it. With this choice of jumps the first dragnet of size $d$ materializes at $j=d^2+1,$ the beginning of the first interval, which is also the first jump site. The first interval ends at $(d+1)^2$. With the Algorithm one can always generate at least sequences of length $d^2$, but after that i.e. at the beginning of the first interval there is the possibility that the sites in the dragnet generate a full block which terminates the generation. 

\vskip .5truein

\noindent {\subtitle 2.2. A probabilistic model}

\vskip .2truein
\noindent We know proceed to set-up and investigate a simplified model for the halting of the Algorithm.

To state our first result we first recall that for non-negative integers $k_i$ the {\bf multinomial coefficient} is defined by the expression
$$\left({k_1+k_2+\cdots +k_d}\atop {k_1\ k_2\ \ldots\ k_d}\right)={(k_1+k_2+\cdots +k_d)!\over {k_1!k_2!\ldots k_d!}}$$
\vskip .1truein
\noindent where $0!=1$ if need be. If $d=2$ this is just the Binomial coefficient. The following basic property of multinomials will be needed subsequently:

\proclaim The Multinomial Theorem: For $n\ge d$
$$\sum_{{k_1,k_2,\ldots,k_d}\atop {k_1+\cdots+ k_d=n}}\left(n\atop {k_1\ k_2\ \ldots\ k_d}\right)=d^n\leqno{(2.1)}$$
where the sum is $d$-fold over the non-negative integers $k_i$.

\noindent The proof amounts to inductively rewriting $(\sum_1^d1)^n$ in a more explicit form. It can be found in many combinatorics and statistics books (see e.g. [A]).

\proclaim Theorem 2.3.: Assume that all the symbols on $\{1,2,\ldots,j-1\}$ have been laid out independently and uniformly from $S$ i.e. they are distributed $\sim B\left({1\over d},{1\over d},\ldots,{1\over d}\right).$ Let $B_j$ be the event that one has the first full block at $j$ in the $i^{\rm th}$ interval. Then
$${\rm {\bf P}}(B_j)=p_i={1\over d^{d+i-1}}\sum_{{k_r\ge 1,\ r=1,\ldots, d}\atop {k_1+\cdots +k_d=d+i-1}}\left({d+i-1}\atop {k_1\ k_2\ \ldots\ k_d}\right)\leqno{(2.2)}$$
where the sum is $d$-fold over the given positive integers.\par

\vskip .1truein
\noindent {\bf Proof:} For $j$ in the $i^{\rm th}$ interval the dragnet is of cardinality $d+i-1.$ Since the entries on the support of the dragnet $\left(j_1,\ldots,j_{d+1-1}\right)$ are independent, each elementary event $\left(s_{j_1},\ldots,s_{j_{d+1-1}}\right)$ is of probability $d^{-(d+i-1)}.$ A full block at $j$ materializes if and only if each symbol of $S$ appears at least once in this dragnet. Since the order of the symbols is immaterial the multiplicity of the elementary event is given by the multinomial coefficient
$$\left({d+i-1}\atop {k_1\ k_2\ \ldots\ k_d}\right)\quad{\rm with}\quad k_r\ge 1\ {\rm for\ all}\ r=1,\ldots,d.$$
All these arrangements on the dragnet are disjoint, hence the sum (2.2) gives the total probability of the full block. \hfill\QED

\vskip .2truein
\noindent {\bf Remark:} The independence assumption here could of course be a heavy one. For one thing identical subsequent symbols ($s_j=s_{j+1}$) are forbidden in the original model but not here. For large $S$ one expects exclusions like this to be less significant i.e. our simplified model should then be more accurate. We return to this in Section 2.2. after the data comparison.  

\vskip .2truein
\noindent The sequence $\{p_i\}$ is crucial in understanding our model. It is of course non-decreasing since the dragnet is not shrinking with increasing $i$. Indeed one can easily show that $p_i$ is strictly increasing and $p_i\uparrow 1.$ More specifically

\proclaim Proposition 2.4.: For an alphabet $S$ of size $d$ one has for all $i\ge 1$
$$1-p_i< d\left(1-{1\over d}\right)^{d-1}\left(1-{1\over d}\right)^i.\leqno{(2.3)}$$ 
\par

\vskip .1truein
\noindent {\bf Proof:} Geometrically the Multinomial Theorem (2.1) is a statement about the total sum on the entries on the $n^{\rm th}$ level of Pascal's $d$-simplex i.e. on the simplex face perpendicular to $(1,1,\ldots,1)\in {\rm {\bf N}}^d$ at (lattice) distance $n$ from the top of the pyramid (which is at the origin). For the size of $1-p_i$ we need by Theorem 2.3 to estimate the boundary sum on that face:
$$\sum_{{k_r=0\ {\rm for\ some\ }r}\atop {k_1+\cdots +k_d=d+i-1}}\left({d+i-1}\atop {k_1\ k_2\ldots\ k_d}\right)$$
Without loss of generality assume that at least $k_d=0.$ Then by the Multinomial Theorem
$$\sum_{{k_1,\ldots,k_{d-1}}\atop {k_1+\cdots +k_{d-1}=d+i-1}}\left({d+i-1}\atop {k_1\ k_2\ldots\ k_{d-1}}\right)=(d-1)^{d+i-1}.\leqno{(2.4)}$$
But $d$-simplex has $d+1$ sides and its $(d+i-1)^{\rm st}$ level has $d$ sides. For each of them (2.4) holds, hence for all boundary terms we have the bound
$$\sum_{{k_r=0\ {\rm for\ some}\ r}\atop {k_1+\cdots +k_d=d+i-1}}\left({d+i-1}\atop {k_1\ k_2\ldots\ k_d}\right)<d(d-1)^{d+i-1}\leqno{(2.5)}$$
which is a strict inequality since corners are accounted multiple times and the multinomial is at least $1.$ But dividing (2.5) by $d^{d+i-1}$, the total sum over the $(d+i-1)^{\rm st}$ level, gives the result. \hfill\QED

\vskip .2truein
\noindent {\bf Remark:} Although this is not subsequently used, it is interesting to note that in spite of the slight slack in (2.5), there seems to be no asymptotic gap in (2.3) (numerical check).  

\vskip .2truein

\proclaim Theorem 2.5.: Let the assumption on the sequence be as in Theorem 2.3. and let $f(n)=n^2.$ Then the sequence generation halts at $j$ i.e. the distribution of the full block is given by
$${\rm {\bf P}}({\rm halts\ at\ } j)=\cases{0,& $1\le j\le d^2$;\cr
(1-p_1)^{j-d^2-1}p_1, & $j$ in the first interval;\cr
\left(\prod_{k=1}^{i-1}(1-p_k)^{l_k}\right)(1-p_i)^{j-d^2-1-\sum_{k=1}^{i-1}l_k}p_i, & $j$ in the $i^{\rm th}$ interval, $i\ge 2$.\cr}\leqno{(2.6)}$$
\par

\vskip .1truein
\noindent {\bf Proof:} Before $j=d^2+1$ the dragnet is too small to support all of $S$ hence full block is impossible. On the first interval the distribution is geometric: with probability $1-p_1$ no block and with probability $p_1$ a block and thereby halting. At each instant the termination is independent of the past.

At times $j=(d+i)^2+1,\ i=2,3,\ldots $ the parameter changes to from $p_i$ to $p_{i+1}$ but otherwise the termination mechanism is still geometric. Compounding this to account the later intervals gives the general case in the expression (2.6).\hfill\QED

\vskip .2truein
\proclaim Corollary 2.6.: The halting time distribution has an exponential tail. The sequences generated are almost surely finite.
\par

\vskip .1truein
\noindent {\bf Proof:} Within the $i^{\rm th}$ interval the halting probability (2.6) is geometrically decreasing with the rate $1-p_i$. Between the $i^{\rm th}$ and  $(i+1)^{\rm st}$ intervals the probability jumps,
$${{{\rm {\bf P}}({\rm halts\ at\ } (d+i)^2+1)}\over {{\rm {\bf P}}({\rm halts\ at\ } (d+i)^2)}}={{(1-p_i)p_{i+1}}\over {p_i}}=r_i\ .\leqno{(2.7)}$$
For small $i$ this expression can and will exceed 1. But by Proposition 2.4. $p_i\rightarrow 1$ and hence also ${p_{i+1}\over p_i}\rightarrow 1$ so $\max{\left\{r_i,1-p_i\right\}}$ will eventually remain below $1-\epsilon$ for some $\epsilon>0$ which implies the result for $i$. The lengths $l_i$ satisfy $l_{i+1}=l_i+2,$ so in particular $l_{i+1}\le (1+\delta)l_i$ for some $\delta>0$ and the exponential decay holds in $j$ as well.

If $A_j=\left\{\omega\ |\ {\rm sequences\ of\ length} \ge j\right\}$ then by the above ${\rm {\bf P}}\left(A_j\right)\le e^{-cj}$ for some positive $c.$ Hence $\sum_j{\rm {\bf P}}\left(A_j\right)$ is summable and by the Borel-Cantelli Lemma ${\rm {\bf P}}\left(A_j\ {\rm i.o.}\right)=0.$
\hfill\QED

\vskip .2truein
\noindent {\bf Remark:} All the results above before Theorem 2.5. were independent of the particular choice $f(n)=n^2.$ In this theorem the assumption was made just to have an explicit (but still complicated) formula. In the Corollary one needs to control the interval lengths well enough, a task difficult only for rather wildly behaving $f.$

\vfill
\eject

\vskip .5truein
\noindent {\subtitle 2.3. Model versus data}
\vskip .2truein

\noindent To investigate the properties of the sequences generated by the Algorithm we collected samples for alphabets of sizes $d=4,5,7,10,15$ and $20.$ The upper limit arises simply from the computational effort. If we set $M=2000,$ a reasonable guess (after some playing around) for an upper bound for termination for $d=20$, the sequences must be sieved out from $20^{2000}$ candidates using the Algorithm. Computed with {\sl Mathematica} on a decent workstation these runs are timed in days. The sample sizes (number of sequences generated) that the following observations are based are indicated in Table I, last column. 

On the other hand to assess the probabilistic model the main hurdle is to compute explicit values for $p_i$ for sufficiently many $i=1,2,3,\ldots$ since these parametrize the halting distribution in (2.6). Here one can expedite things by computing the multinomials in the most efficient way. The key sum in (2.2),

$$\sum_{{k_r\ge 1,\ r=1,\ldots, d}\atop {k_1+\cdots +k_d=d+i-1}}\left({d+i-1}\atop {k_1\ k_2\ldots\ k_d}\right),$$

\noindent can be simplified. E.g. for $i=4$ it can be computed as
$$\left({d\atop 1}\right)\left({d+3}\atop {1\ \ldots\ 1\ 4}\right)+2\left({d\atop 2}\right)\left({d+3}\atop {1\ \ldots\ 1\ 2\ 3}\right)+\left({d\atop 3}\right)\left({d+3}\atop {1\ \ldots\ 1\ 2\ 2\ 2}\right),$$
an expression far easier to deal with than the original ( $\ldots\ $ denotes here a constant string of 1's). The idea here is to keep record of the multiplicity of the string that appears inside the multinomial (like $23$ in the middle expression - it can appear in $\left({d\atop 2}\right)$ different ways and in 2 different orders). But this too gets rather unwieldy for large $d$ with $i$ beyond the smallest values.

\vskip .2truein
\noindent Figure 1., middle row, is joint plot of the sequence termination instant distributions of the data and the model for $d=10$. One sees a reasonable overall match but also distinct discrepancies due to our simplified model. Both plots peak at the same spot, the start of the fifth interval ($i=5\ {\rm so\ } (d+i-1)^2+1=197$) and their levels there are fairly close to each other (0.0140179 and 0.0150393, empirical and model resp.). The model favours slightly earlier halting as witnessed in the distributions: it both rises and sets earlier.

The model graph has by definition geometric slopes on the intervals and so indeed does the data: this is clear from the log-plot on the right of the same data (plot is $\ln({\rm number\ of\ samples} + 1)$). The experimental and model slope tilts are actually fairly close to each other as seen on the left, especially in the middle part of the distribution.  
The key observation here though is that the tail of the data distribution seems indeed majored by an exponential and tapers off around same range of $j$ as it does for the model. This strongly suggests that the sequences don't get very far. 
The top row of Figure 1 present the same comparison for $d=5.$ The overall distribution match is actually quite remarkable. In the finer details like the slope roughness there is deviation from the model (although the sample size is here at its largest, see Table 1). We return to this in Section 2.3. 

\vfill
\eject

\vskip .8truein
\vbox{
\hbox{
 \psfig{figure=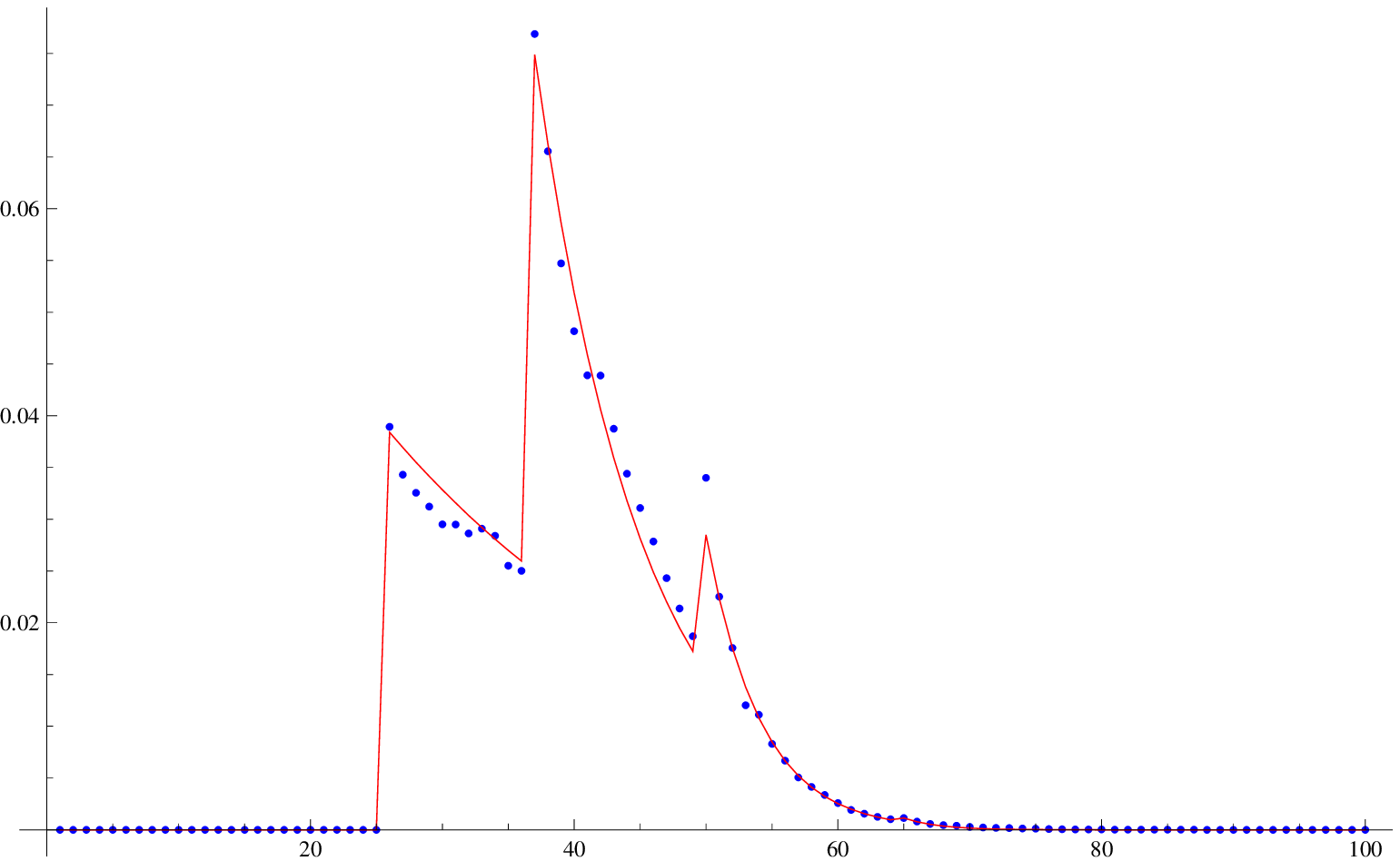,height=1.5in}
 \hskip .3in
 \psfig{figure=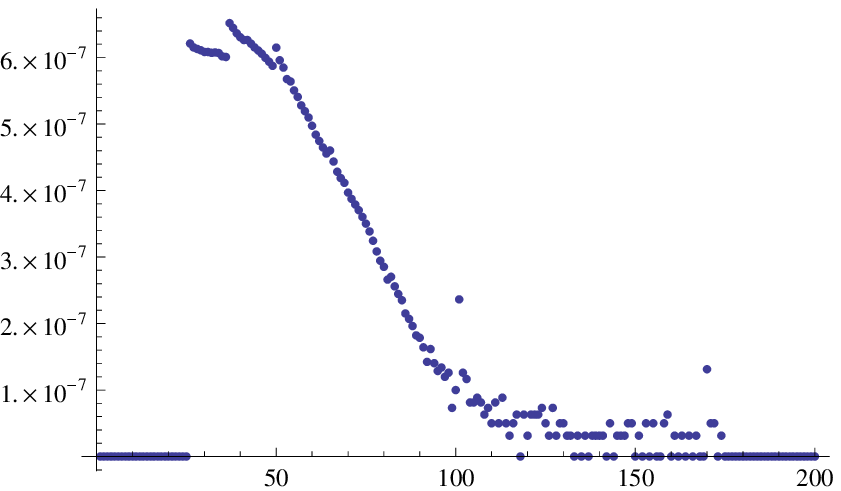,height=1.5in}
}
\vskip .3truein
\hbox{
 \psfig{figure=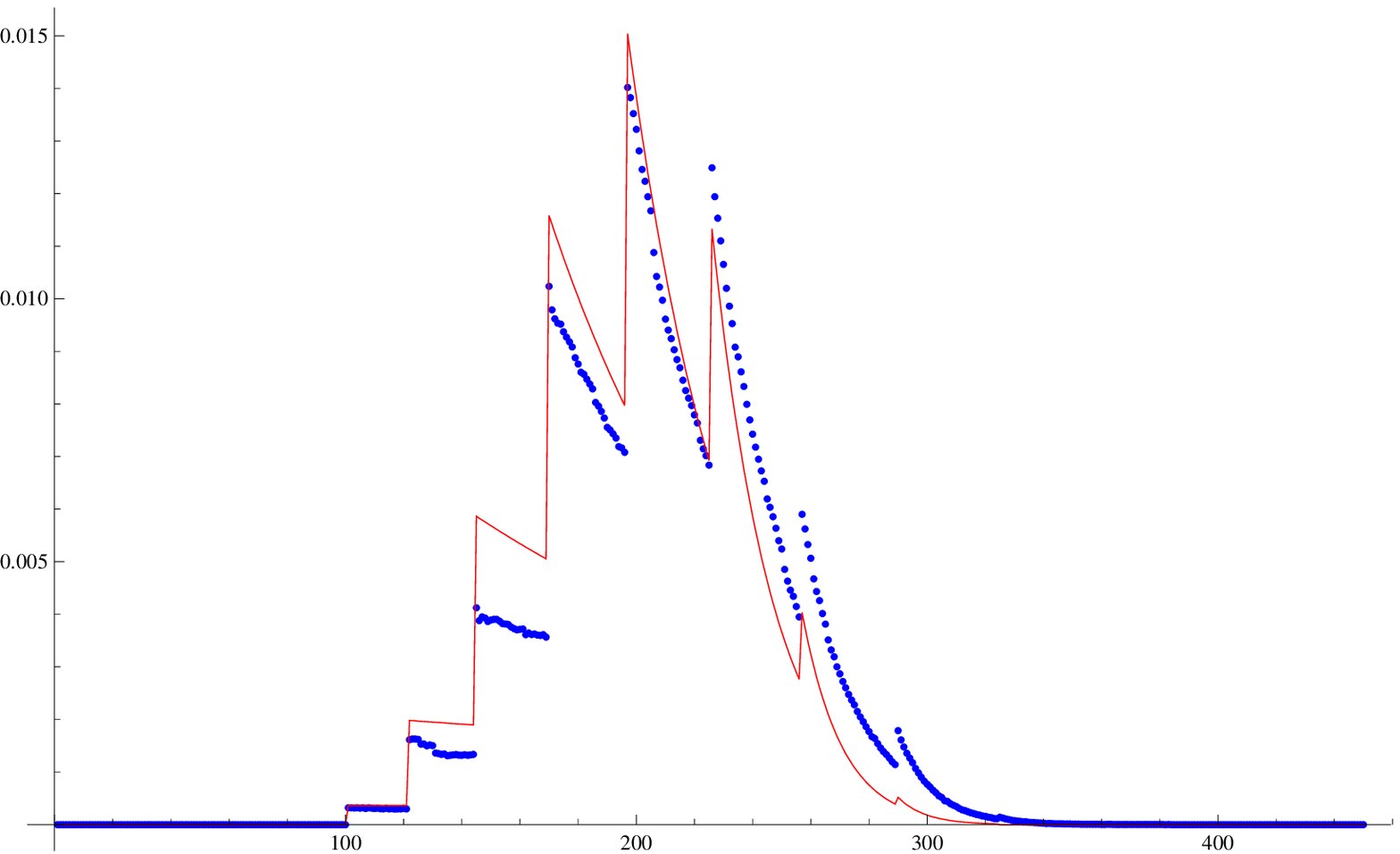,height=1.5in}
 \hskip .3in
 \psfig{figure=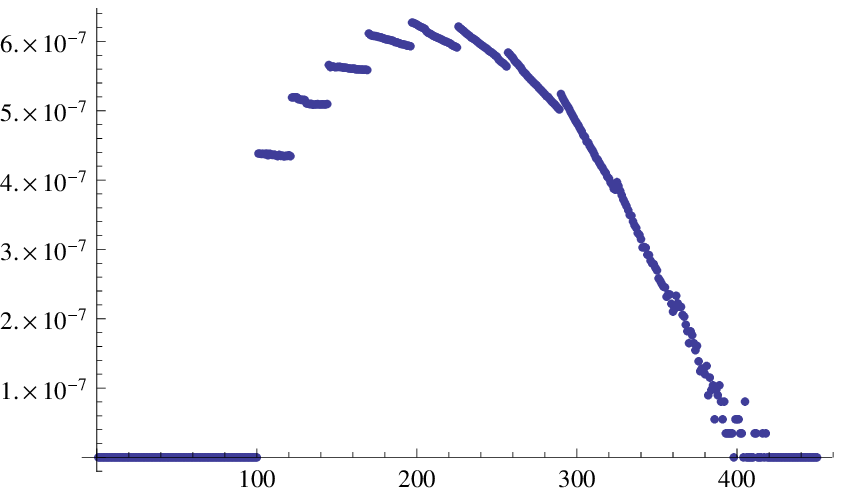,height=1.5in}
}
\vskip .3truein
\hbox{
 \psfig{figure=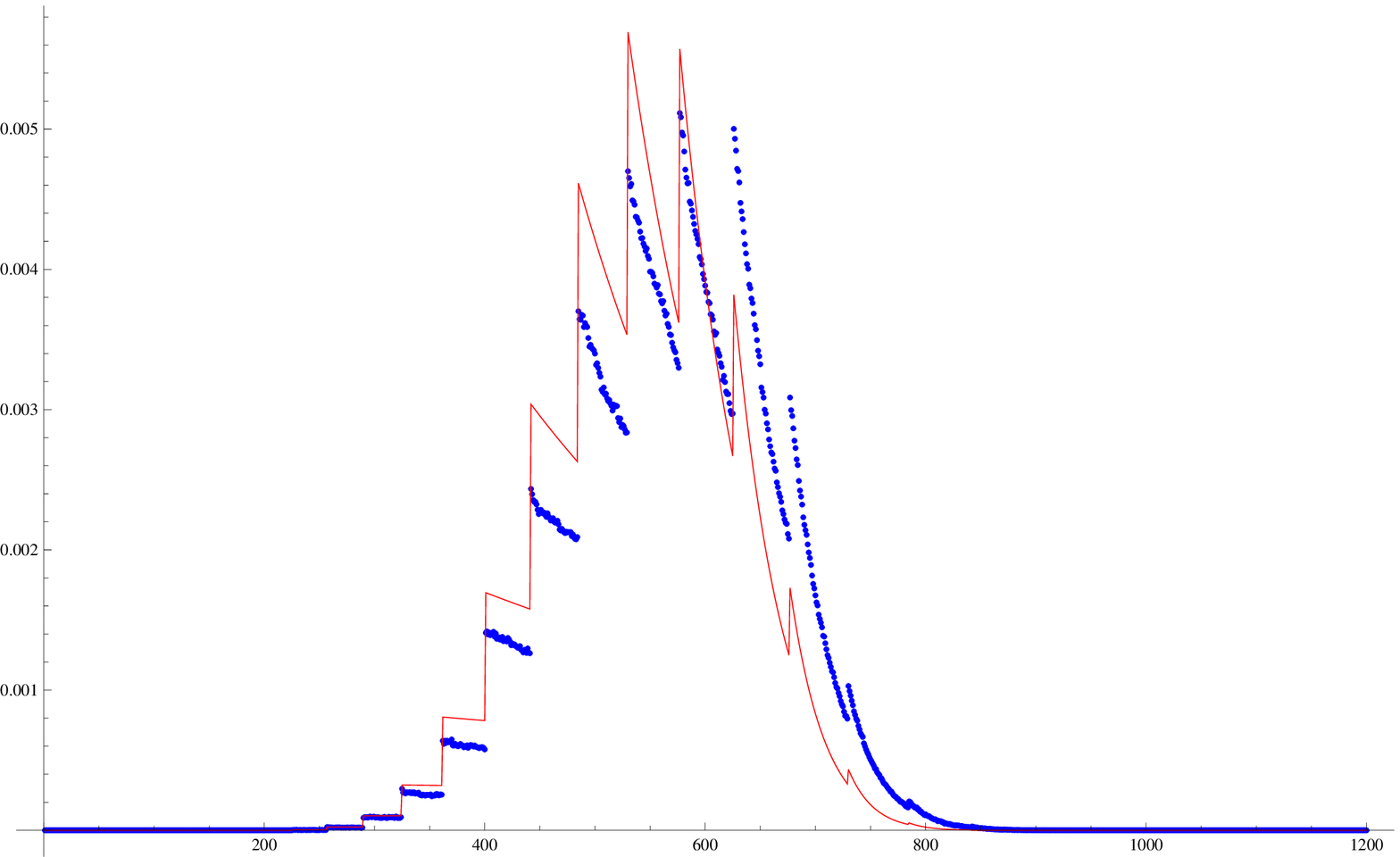,height=1.6in}
 \hskip .3in
 \psfig{figure=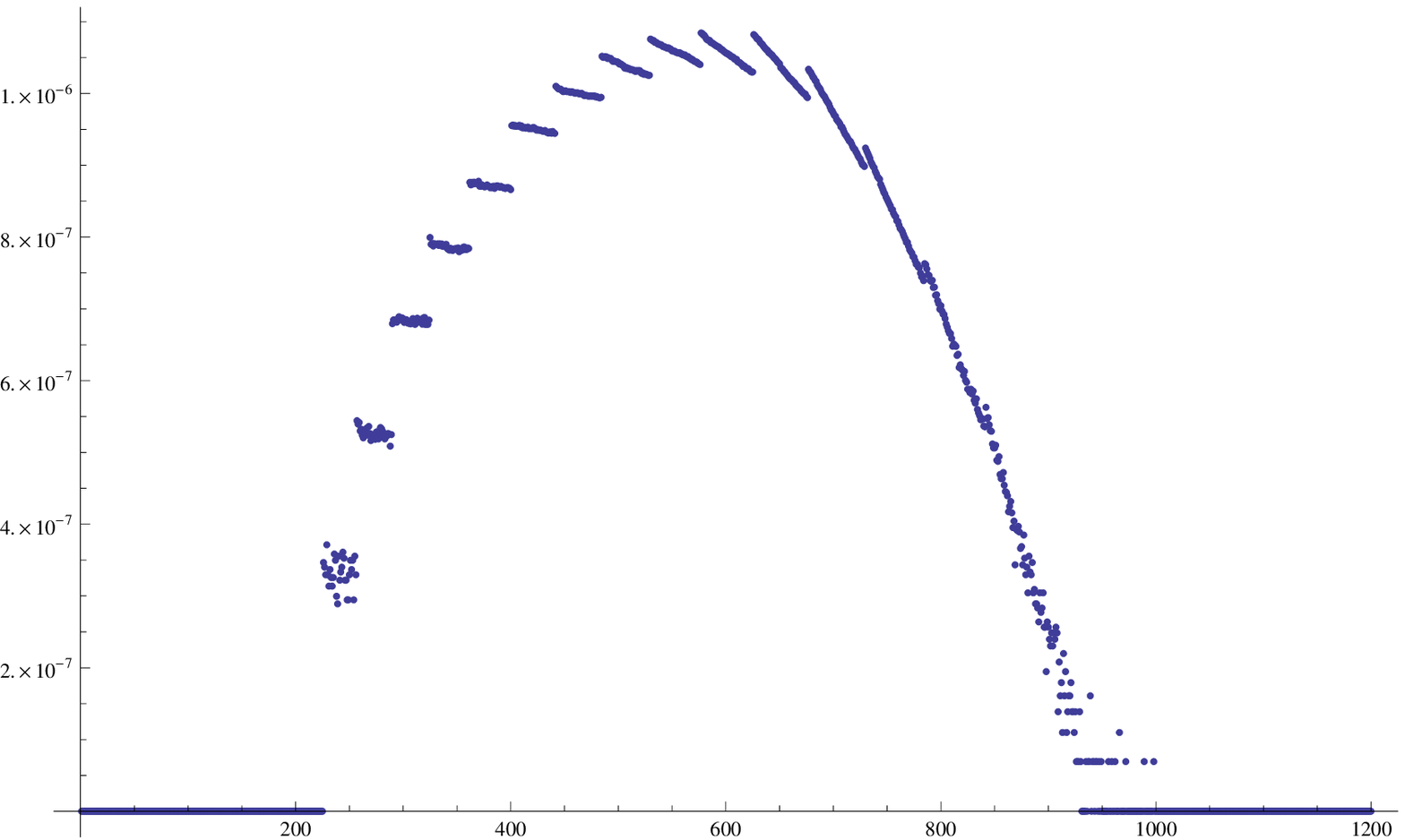,height=1.6in}
}
}

\vskip .2truein
\noindent{\bf Figure 1:} Empirical (blue/rough) and theoretical (red/continuous) halting distributions superimposed together (left) and corresponding log empirical distribution (right) for different alphabet sizes. From top to bottom row: $d=5,10$ and $15.$

\vskip .2truein
\noindent The left plots of the bottom row (for $d=15$) exhibits increasing symmetry in the distributions. It is also hinting at the high $d$-limit: with suitable scaling the halting distribution might be normal. Magnifying the slopes here shows that the mechanism that gave the top row data the \lq\lq roughness\rq\rq\ are not there anymore. One is tempted to think that this is due to randomness smoothing them out.

\vskip .2truein
\noindent Perhaps a bit troubling facts in the mid and bottom plots on the left in Figure 1 are the vertical off-set between the slopes. There is a way to see a dependency mechanism in the sequences that contributes to them appearing. While this does not seem to provide a way to improve the model it is still worth understanding.

Consider the assignment sequence indicated in Figure 2. Here $a,b,c$ and $d$ are coordinates, $b=a+k^2,\ c=a+n^2,\ d=b+n^2,\ n>k.$ Since $x_a=s$ we know that $x_b\not=s\not=x_c$ as indicated. One should think of a dragnet being at $d$, hence $b$ and $c$ belong to it and $a$ may or may not belong. The simplest case (and probably the one contributing most to the probabilities involved) is of course $k=1,\ n=2.$ We are interested in the probability of a full block appearing at $d.$

\vskip.4truein
\centerline{\hbox{
 \psfig{figure=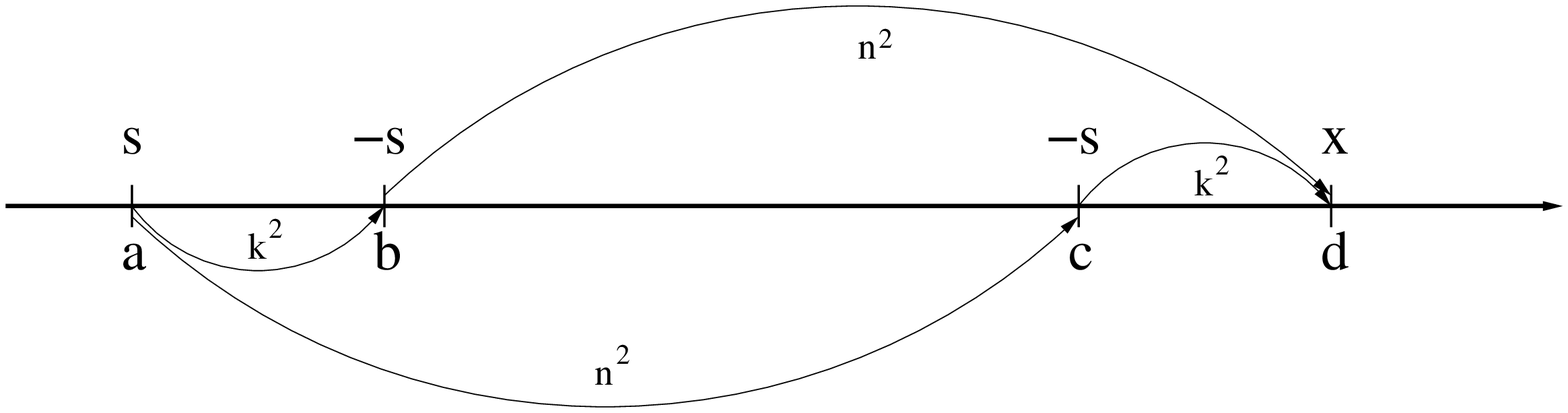,height=1in}
}}

\vskip .1truein
\noindent{\bf Figure 2:} A dependency mechanism affecting the termination probability.

\vskip .2truein
\noindent Ignoring all other symbols on the interval $\{a,\ldots,d\}$ two cases should be distinguished:
\vskip.1truein
\item{1.} if $k^2+n^2$ is not a square then no symbol from the triple $x_a,x_b\ {\rm and\ }x_c$ forbids $x_d=s$ but $x_b$ and $x_c$ do block one or two non-$s$ symbols at $d.$
\item{2.} if $k^2+n^2$ is a square then $x_a$ forbids $x_d=s$ and furthermore $x_b$ and $x_c$ block one or two  non-$s$ symbols at $d.$

\vskip.1truein
\noindent On the basis of these we see immediately that
$${\rm {\bf P}}\left({\rm full\ block\ at\ } d\ |\ k^2+n^2\ {\rm not\ square} \right)< {\rm {\bf P}}\left({\rm full\ block\ at\ } d\right)< {\rm {\bf P}}\left({\rm full\ block\ at\ } d\ |\ k^2+n^2\ {\rm square} \right).$$
\noindent But as we vary $k$ and $n$ there are many more sites $a$ where 1. instead of 2. holds. Hence one should expect the full block i.e. termination probabilities be majored by the independent probabilities of the model. This is indeed the case as seen in Figure 1. How to quantify this is another matter though, well worth further study.

\vskip .4truein
\noindent To further characterize the sampled sequences we computed their expected lengths and standard deviations. These and the corresponding model data are presented in Table I for a few $d$ values. For large $d$ we already mentioned of the problem of computing $p_i$ for high $i$. These are needed for an accurate tail estimate of the distribution. The asterisks in the Table are due to this complication.

The values of the two statistical indicators here just quantify aspects of the earlier distributional matching. A bit surprisingly the best match occurs at the value $d=5$ indicating that small alphabets don't seem to suffer from the independence assumption in the model.

The decay rate of the density in Theorem 1.9. can be used to estimate termination times. These turn out truly staggering compared to the real ones above. This is of course due to the proof methods and suggests that other arguments truer to the nature of the problem should be found.  

\hfill
\eject

\vskip .5truein
{
\offinterlineskip
\tabskip=0pt
\halign{ 
\vrule height2.75ex depth1.25ex width 0.6pt #\tabskip=1em
& #\hfil & \vrule # \hfil &  #\hfil & \vrule # & # \hfil & # \hfil \vrule & # \hfil & \vrule # & #\hfil & \vrule # & #\hfil & #\vrule width 0.6pt \tabskip=0pt\cr
\noalign{\hrule height 1pt}
& Symbols && Empirical && Empirical && Model && Model && Sequences &\cr
& $d$ && mean && std. dev.  && mean && std. dev. &&  &\cr
\noalign{\hrule height 1pt}
& 4 && 27.2542 && 5.13374 && 23.992 && 5.23924 && $50\cdot 10^6$ &\cr
\noalign{\hrule}
& 5 && 39.5672 && 8.28983 && 39.2172 && 8.22516 && $80\cdot 10^6$ &\cr
\noalign{\hrule}
& 6 && 60.8247 && 13.5813 && 59.3666 && 11.9713 && $80\cdot 10^6$ &\cr
\noalign{\hrule}
& 7 && 89.4687 && 18.5912 && 84.982 && 16.5113 && $30\cdot 10^6$ &\cr
\noalign{\hrule}
& 10 && 209.315 && 38.2887 && 199.562 && 35.1369 && $20\cdot 10^6$ &\cr
\noalign{\hrule}
& 15 && 566.87 && 92.2796 && 543.291 && 84.4349 && $10\cdot 10^6$ &\cr
\noalign{\hrule}
& 20 && 1156.57 && 170.829 && $\ast$ && $\ast$ && $5\cdot 10^6$ &\cr
\noalign{\hrule height 1pt}
}}

\vskip .3truein
\noindent {\bf Table I}. Some empirical and theoretical statistics of the one-sided sequences. Means and standard\hfill\break deviations of the halting times at various alphabet sizes. Asterisks refer to missing coefficients.

\vskip .2truein
\noindent Analyzing the empirical data further reveals rather clearly certain exponents at play: the values in the second column grow almost exactly like $d^{5/2}$ and those of the third column roughly like $d^{15/7}.$ Since there seems to be an increasingly random behavior in the halting of the sequences for higher $d$, we venture to

\proclaim Conjecture 2.7.: Suppose $T^{(d)}$ is the halting instant of the Algorithm v2.0 with $f(n)=n^2.$ For sufficiently rapidly growing $M(d)$ there are positive constants $a$ and $b$ such that as $d\rightarrow\infty$
$${\rm {\bf P}}\left({{T^{(d)}-ad^{5/2}}\over {bd^{15/7}}}\le x\right)\longrightarrow\Phi(x)\qquad \forall x\in {\rm {\bf R}}$$
where $\Phi$ is the cumulative distribution function of the standard normal $N(0,1).$
\par

\vskip .2truein
\noindent {\bf Remarks: 1.} For the sequence length $M(d)$ we just need a rate that outgrows the off-set rate $d^{5/2}.$ Clearly for the left threshold $(d^2+1-ad^{5/2})/bd^{15/
7}\rightarrow-\infty$ as desired. \hfill\break
\noindent {\bf 2.} A partial result would be to show the Central Limit Theorem for the model. For this one would need to analyze the asymptotics of $p_i(d)$ and then the scaling limit of (2.6).
\hfill\break
\noindent {\bf 3.} The CLT for the algorithm would of course imply the a.s. emptiness of the sequence spaces.

\vskip .5truein
\noindent {\subtitle 2.4. Some fine detail}
\vskip .2truein

\noindent While generating the sequences some additional information was gleaned, too. To get some insight in how the haltings come about we recorded the pairs $(i,n)$ at which the first terminal block was generated. With this one knows that while the sequence may be extended a bit (at least one step) it cannot be continued past $i+n^2.$ Note that because of this the resulting distribution for $i$ is not exactly the distribution of the halting instant $j$ earlier, but slightly shifted to the left.

\vskip .2truein
\noindent {\bf Algorithm v2.1:}\hfill\break 
\indent {\tt 0. set} $M\ge 1,$ {\tt let} $S_j=S$ {\tt at each}  $j\in \{1,\ldots,M\}$ {\tt and set} $i=1.$\hfill\break
\indent {\tt 1. pick uniformly a random symbol} $s\in S_i.$ \hfill\break
\indent {\tt 2. update} $S_j\leftarrow S_j\setminus \{s\}$ {\tt for\ all} $j=i+f(n)\in \{i+1,\ldots,M\},\ n\in {\rm {\bf N}}$ \hfill\break
\indent\hskip.75cm {\tt if} $S_j=\emptyset$ {\tt results for some $j$, record} $(i,n),$ {\tt halt\ and\ call} {\bf empty},\hfill\break
\indent\hskip .75cm {\tt if} $i=M$ {\tt halt\ and\ call} {\bf full\ length},\hfill\break
\indent\hskip.75cm {\tt else} $i\leftarrow i+1$ {\tt and go to 1.}

\vskip .2truein
\noindent  With this modified algorithm one gets to see from where the terminal jumps originate. Figure 3. illustrates their distributions, call it the $(i,n)${\bf -plot}, for two $d$-values.

\vskip .5truein
\vbox{
\hbox{
 \psfig{figure=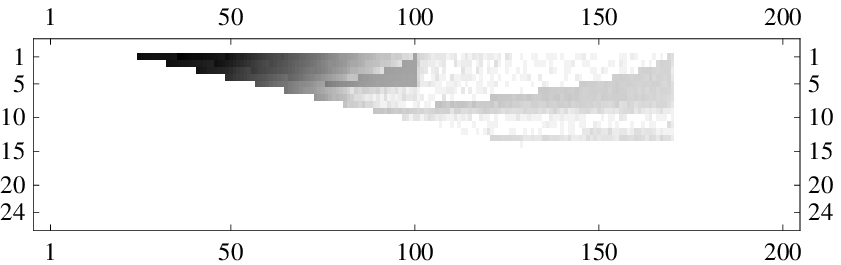,height=.88in}
 \hskip .3in
 \psfig{figure=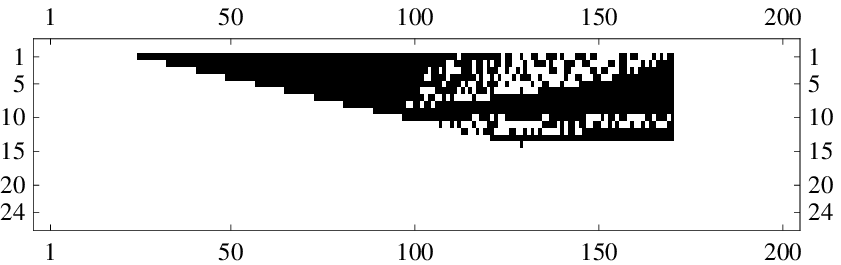,height=.88in}
}
\vskip .3truein
\hbox{
 \psfig{figure=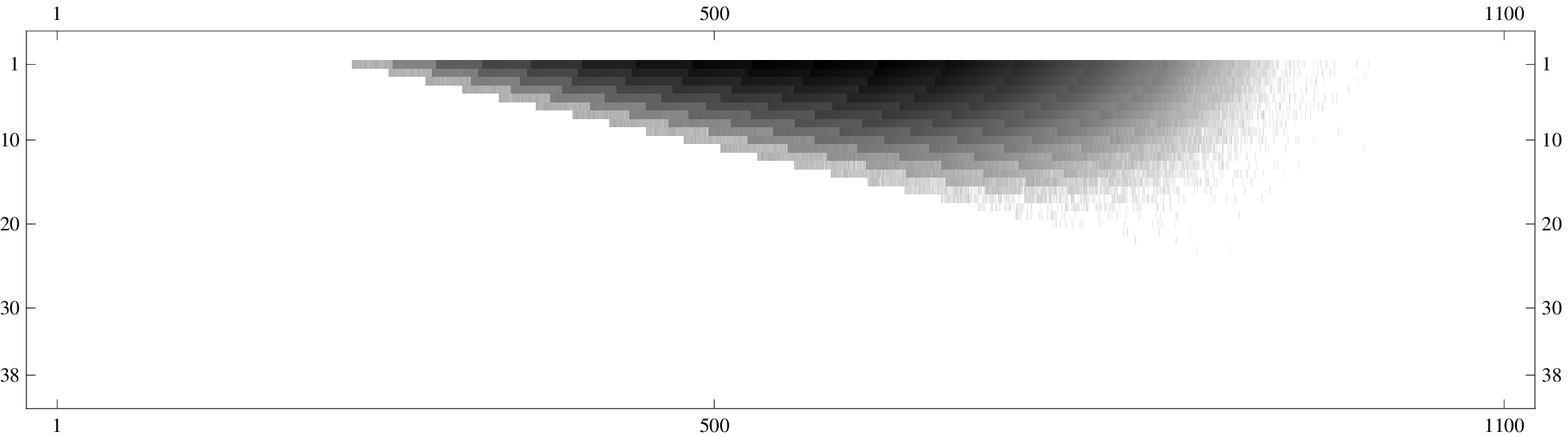,height=.77in}
 \hskip .3in
 \psfig{figure=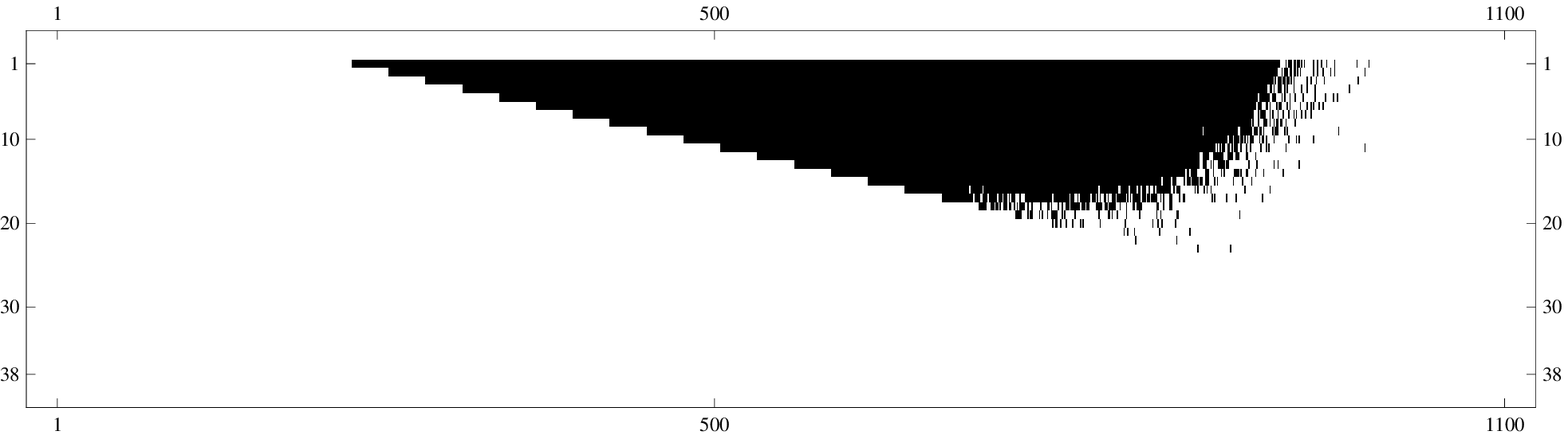,height=.77in}
}
}

\vskip .2truein
\noindent{\bf Figure 3:} Termination distributions for $d=5$ (top) and $15$ (bottom). Iterate $i$ runs rightwards and jump size $n$ downwards. Same distributions in each row; on the left a log-plot, on the right the plain support (non-zero entries). Data from $20$ and $10$ million samples respectively. Note the different scales up and down.

\vskip .2truein
\noindent The first distinct feature, the left hand \lq\lq staircase\rq\rq\ in the plots, can be readily explained.

\proclaim Proposition 2.8.: The staircase jumps in the $(i,n)$-plot are at $d^2+2(n-1)d-2(n-1),\ n\ge 1$ and the step length is $2(d-1).$  
\par

\vskip .1truein
\noindent {\bf Proof:} If a terminal block is generated with a jump of size $n^2$, then the previous $d-1$ (non terminal) blocks into this future site have been generated with jumps of size at least $(n+1)^2,(n+2)^2,\ldots ,(n+d-1)^2.$ These particular jump sizes give the minimal $i$ for which $i+n^2$ results in a termination. So necessarily $j=i+n^2\ge (n+d-1)^2+1,$ from which we get that $i\ge (n+d-1)^2-n^2+1=\cdots =d^2+2\left((n-1)d-n+1\right).$ Evaluating this for $n\ge 1$ gives the step values and their difference gives the step size. \hfill\QED

\vskip .2truein
\noindent In a similar fashion to Definition 2.1. one can formulate locations and length of the horizontal $(l\times 1)$-blocks inside the distributions at any given level $n.$ Within these the distribution of terminations seems to be approximately geometric as expected from our earlier analysis.

A couple of things in these plots are far more intriguing though. The first is the internal and right end structure in the plots for $d=5.$ There seem to be some hard combinatorial constraints at play here. Hence the lace-like interior and the abrupt wall on the right (at 173). Note that the plots are based on a large data set of $20\cdot 10^6$ samples. Furthermore, for the computation the sequence length attempted was $M= 600$ i.e. the shown graphics is a cutout without any boundary effect from 200. Also a blank i.e. no sample points is rather improbable event for low $n$ in the interior of the plot unless there is a definite reason for the absence. 

The lower plots clearly indicate random phenomena dominating for higher a alphabet size (here $d=15$). Apart from the left staircase there are no traces of hard constraints elsewhere. The right hand side of the distribution is tapering off like the halting distribution earlier, most likely with exponential tails at all $n$-levels. 


\vskip .5truein
\noindent {\subtitle 3. Conclusions}
\vskip .2truein

\noindent Here we have tried to chart of what still seems unexplored or at least rather hazy terrain. This of course applies to much of the study on long range order e.g. in Statistical Mechanics and related fields. To come up with an empty set is perhaps a bit disappointing but hopefully the methods here are useful in further work. At least we know how stringent even the most innocuous looking exclusion rules are if applied with unbounded reach (once the obvious case of periodicity has been ruled out). Moreover the large discrepancy between the termination times arising from proofs and the simple probabilistic model suggests that there should be a more natural proof to be found.

In the latter half we have restricted ourselves mostly for the jump sequence $f(n)=n^2.$ The reason is two-fold: firstly, the original question was about this and secondly other nontrivial $f$'s get quickly out of hand computationally. Our probabilistic model should work verbatim for $n^k,\ k\ge 3$ but because of their jump sizes, computing the exact values of $p_i$ for high $i$ remains a problem. So does generating data: the range $M$ in our algorithms has to be much larger hence a far bigger number of exclusions have to be recorded. The reason why at high $d$ a sequence generated eventually terminates remains the same, the bulk accumulation of exclusions at a site (which is overcome only beyond the lacunar threshold). But a much larger store of excluded symbol values need to be kept in order to finally record a full block. Conjecture 2.7. is likely to be true more generally but because of lack of supporting statistical data we will not venture to guess it.

\vskip .5truein
\noindent {\subtitle Acknowledgements}
\vskip .2truein

\noindent The author would like to thank Mike Keane, Rob Tijdeman and Benjamin Weiss for valuable comments and pointers that led to the right literature. 

\vfill
\eject

\vskip .4truein
{\subtitle References}

\vskip .2truein

\item{[A]} Aigner, M.: {\sl A course in Enumeration}, Springer, 2007.
\item{[AHK]} Ajtai, M., Havas, I, Koml\'os, J.: Every group admits a bad topology, {\sl Studies in Pure Mathematics}, pp. 21-34, Birkh\"auser, 1983.
\item{[BPPS]} Balog, A., Pelik\'an, J., Pintz, J., Szemer\'edi, E.: Difference sets without k-th powers, {\sl Acta Math. Hungar.}, {\bf 65} (2), pp. 165-187, 1994
\item{[B-M]} Betrand-Mathis, A.: Ensembles intersectifs et recurrence de Poincar\'e, {\sl Colloque de Th\'eorie Analytique des Nombres \lq\lq Jean Coquet\rq\rq (Marseille, 1985)}, pp. 55-72, Publ. Math. Orsay, 1988-2, Univ. Paris XI, Orsay, 1988.
\item{[CT]} Chow, Y.S., Teicher, H.: {\sl Probability Theory: Independence, Interchangeability, Martingales}, Springer, 1978.
\item{[F]} Furstenberg, H.: Ergodic behavior of diagonal measures and a theorem of Szemer\'edi on arithmetic progressions, {\sl J. d'Analyse Math.}, {\bf 71}, pp. 204-256, 1977.
\item{[HU]} Hopcroft, J.E., Ullman, J.D.: {\sl Introduction to Automata Theory, Languages, and Computation}, Addison-Wesley Publishing, Reading Massachusetts, 1979.
\item{[KMF]} Kamae, T., Mend\'es France, M.: van der Corput difference theorem, {\sl Israel J. Math.}, {\bf 31} (3-4), pp. 335-342, 1978.
\item{[K]} Katznelson, Y.: Chromatic numbers of Cayley graphs on ${\rm {\bf Z}}$ and recurrence, {\sl Combinatorica}, {\bf 21}, pp. 211-9, 2001.
\item{[L]} Lucier, J.: Intersective sets given by a polynomial, {\sl Acta Arith.}, {\bf 123}, pp. 57-95, 2006.
\item{[LM]} Lind, D., Marcus, B.:{\sl An introduction to Symbolic Dynamics and Coding}, Cambridge Univ. Press, 1995.
\item{[PS]} Peres, Y., Schlag, W.: Two Erd\"os problems on lacunary sequences: chromatic number and diophantine approximation, {\tt arxiv:math.CO:0706.0223v1}  
\item{[S]} S\'ark\"ozy, A.: On difference sets of sequences of integers, {\sl Acta Math. Hungar.}, {\bf 31} (1-2), pp. 125-149, 1978.
\item{[W]} Weiss, B.: {\sl Single orbit dynamics}, CBMS Regional Conference Series in Mathematics, {\bf 95}, AMS, 2000.
\vfill
\eject

\end